\documentclass{article}
\usepackage{amsfonts}
\usepackage{amsmath}
\usepackage{amssymb}
\usepackage{graphicx}
\numberwithin{equation}{section}
\newtheorem{Theorem}{Theorem}[section]
\newtheorem{Proposition}{Proposition}[section]
\newcommand{\gata}{\blacksquare}
\newtheorem{Remark}{Remark}[section]

\begin{document}

\title{On the geometry of Siegel-Jacobi domains}
\author{S. Berceanu, A. Gheorghe \\
National Institute for Nuclear Physics and Engineering,\\
P.O.Box MG-6, RO-077125 Bucharest-Magurele, Romania\\
E-mail: Berceanu@theory.nipne.ro}
\maketitle

\begin{abstract}
We study the holomorphic unitary representations of the Jacobi group based
on Siegel-Jacobi domains. Explicit polynomial orthonormal bases of the Fock
spaces based on the Siegel-Jacobi disk are obtained. The scalar holomorphic
discrete series of the Jacobi group for the Siegel-Jacobi disk is
constructed and polynomial orthonormal bases of the representation spaces
are given.

\textbf{M.S.C. 2000:}  22E30, 20G05,  11F50, 12E10, 81R30

\end{abstract}
\section{Introduction}\label{sec1}

The Jacobi groups are semidirect products of appropriate semisimple real
algebraic group of Hermitian type with Heisenberg groups \cite{TA99}, 
\cite{LEE03}. The semisimple groups are associated to Hermitian symmetric
domains that are mapped into a Siegel upper half space by equivariant
holomorphic maps \cite{SA80}.  The Jacobi groups are unimodular,
nonreductive, algebraic groups of Harish-Chandra type. The Siegel-Jacobi
domains are nonreductive symmetric domains associated to the Jacobi groups
by the generalized Harish-Chandra embedding \cite{SA80}, \cite{LEE03}, 
\cite{Y02} -\cite{Y08}.

The holomorphic irreducible unitary representations of the Jacobi groups
based on Siegel-Jacobi domains have been constructed by Berndt, 
B\"{o}cherer, Schmidt, and Takase \cite{BS98}, \cite{BB90}, \cite{TA90}-\cite{TA99}
with relevant topics: Jacobi forms, automorphic forms, spherical functions,
theta functions, Hecke operators, and Kuga fiber varieties.

Some coherent state systems based on Siegel-Jacobi domains have been
investigated in the framework of quantum mechanics, geometric quantization,
dequantization, quantum optics, nuclear structure, and signal processing 
\cite{KRSAR82}, \cite{Q90}, \cite{SH03}, \cite{BERC06}-\cite{BERC08B}.

This paper is organized as follows. In Section \ref{sec2} we present  explicit
formulas for the canonical automorphy factors and kernel functions of the
Jacobi groups and corresponding Siegel-Jacobi domains. In Section \ref{sec3} we
introduce a Fock space of holomorphic functions on the Siegel-Jacobi disk.
We obtain explicit polynomial orthonormal bases for this space and the Fock
spaces with inner products associated to points on the Siegel disk
(Proposition \ref{Pr3.1} and Proposition \ref{Pr3.2}). In Section \ref{sec4} we construct the scalar
holomorphic discrete series of the Jacobi group for the Siegel-Jacobi disk
(Proposition \ref{Pr4.3}). We give polynomial orthonormal bases of the
representation spaces (Proposition \ref{Pr4.5}). Finally, we discuss the link
between the coherent state systems based on Siegel-Jacobi domains and the
explicit kernel functions of representation spaces for Jacobi groups.

\textbf{Notation.} We denote by $\mathbb{R}$, $\mathbb{C}$, $\mathbb{Z}$,
and $\mathbb{N}$ the field of real numbers, the field of complex numbers,
the ring of integers, and the set of non-negative integers, respectively. 
$M_{mn}(\mathbb{F})\approxeq\mathbb{F}^{mn}$ denotes the set of all $m\times
n $ matrices with entries in the field $\mathbb{F}$. $M_{1n}(\mathbb{F})$ is
identified with $\mathbb{F}^n$. Set $M_{n}(\mathbb{F})=M_{nn}(\mathbb{F})$.
For any $A\in M_{mn}(\mathbb{F})$, $^{t}A$ denotes the transpose matrix of
$A$. For $A\in M_{mn}(\mathbb{C})$, $\!\bar{A}$ denotes the conjugate matrix
of $A$ and $A^{\dag}=\,^{t}\bar{A}$. For $A\in M_n(\mathbb{C})$, the
inequality $A>0$ means that $A$ is positive definite. The identity matrix of
degree $n$ is denoted by $I_{n}$. Let $\mathcal{O}(\mathfrak{D,}W\mathfrak{)}
$ denote the space of all $W$-valued holomorphic functions on the connected
complex manifold $\mathfrak{D}$ equipped with the topology of uniform
convergence on compact sets. Here $W$ is a finite dimensional Hilbert space.
Set $\mathcal{O}(\mathfrak{D)=}\mathcal{O}(\mathfrak{D,}W\mathfrak{)}$ for
 $\dim W=1$. In this paper we will use the words " unitary representation on a
Hilbert space" to mean a continuous unitary representation on a complex
separable Hilbert space.

\section{Canonical automorphy factor and kernel \\ function for Jacobi groups}\label{sec2}

We begin with the definition of the Jacobi group given in \cite{BS98}, 
\cite{TA99}, \cite{LEE03}. Let $\mathfrak{H}_{n}$ be the Siegel upper half space
of degree $n$ consisting of all symmetric matrices $\Omega 
\in M_{n}(\mathbb{C)}$ with $\operatorname{Im}\Omega >0$. Let $\mathrm{Sp}(n,\mathbb{R})$ be the
symplectic group of degree $n$ consisting of all matrices $\sigma \in 
M_{2n}(\mathbb{R)}$ such that $^{t}\sigma {J_{n}}\sigma =J_{n}$, where 
\begin{equation}
\sigma =\left( 
\begin{matrix}
a & b \\ 
c & d%
\end{matrix}%
\right) ,\quad J_{n}=\left( 
\begin{matrix}
0 & I_{n} \\ 
-I_{n} & 0%
\end{matrix}%
\right) ,  \label{a1}
\end{equation}%
and $a$, $b$, $c$, $d\in M_{n}(\mathbb{R)}$. The group $\mathrm{Sp}(n,\mathbb{
R})$ acts transitively on $\mathfrak{H}_{n}$ by $\sigma \Omega =(a\Omega
+b)(c\Omega +d)^{-1}$, where $\sigma \in \mathrm{Sp}(n,\mathbb{R})$ and $
\,\Omega \in \mathfrak{H}_{n}$. \smallskip

Let $G^{s}$ be a Zariski connected semisimple real algebraic group of
Hermitian type. Let $\mathfrak{D=}G^{s}/K^{s}$ be the associated Hermitian
symmetric domain, where $K^{s}$ is a maximal compact subgroup of $G$.
Suppose there exist a homomorphism $\rho:G^{s}\rightarrow\mathrm{Sp}(n,
\mathbb{R})$ and a holomorphic map $\tau:\mathfrak{D}\rightarrow\mathfrak{H}
_{n}$ such that $\tau(gz)=\rho(g)\tau(z)$ for all $g\in G^{s}$ and $z\in
\mathfrak{D}$. The \emph{Jacobi group} $G^{J}$ is the semidirect product of 
$G^{s}$ and the Heisenberg group $H[V]$ associated with the symplectic $
\mathbb{R}$-space $V$ and the nondegenerate alternating bilinear form $D:$
 $V\times V\rightarrow \mathcal{A}$, where $\mathcal{A}$ is the center of
 $H[V]$. The multiplication operation of $G^{J}\thickapprox G^{s}\times
 V\times\mathcal{A}$ is defined by
\begin{equation}
gg^{\prime}=(\sigma\sigma^{\prime},\rho(\sigma)v^{\prime}+v,\varkappa
+\varkappa^{\prime}+\frac{1}{2}D(v,\rho(\sigma)v^{\prime}),  \label{a2}
\end{equation}
where $g=(\sigma,v,\varkappa)\in G^{J}$, $g^{\prime}=
(\sigma^{\prime},v^{\prime},\varkappa^{\prime})\in G^{J}$.

The \emph{Jacobi-Siegel} \emph{domain }associated to\emph{\ }the Jacobi
group $G^{J}$ is defined by $\mathfrak{D}^{J}=
\mathfrak{D\times}\mathbb{C}^{N}\mathfrak{\cong}G^{J}/(K^{s}\times\mathcal{A})$, where $\dim$ $V=2N$.

Let $w_{0}$ be a fixed element of $\mathfrak{D}$ and let $I_{\tau(w_{0})}$
be the complex structure on $V$ corresponding to 
$\tau(w_{0})\in\mathfrak{H}_{n}$. Let $V_{\mathbb{C}}=V_{+}\oplus V_{-}$ be the complexification of $V$, where $V_{\pm}$ consists of all $v\in V_{\mathbb{C}}$ such that $I_{\tau(w_{0})}v=\pm\mathrm{i}v$. Then $w\in\mathfrak{D}$ and 
$v\in V_{\mathbb{C}}$
 determine the element $v_{w}=v_{+}-\tau(w)v_{-}$ of $V_{+}$, where 
$v=v_{+}+v_{-}$, $v_{\pm}\in V_{\pm}$.

$G^{J}$ is an algebraic group of Harish-Chandra type \cite{TA99},
 \cite{LEE03}, \cite{SA80}. We recall the definition of Harish-Chandra type groups 
\cite{SA80}.

Let $G$ be a Zariski connected $\mathbb{R}$-group with Lie algebra
 $\mathfrak{g}$ and let $G_{\mathbb{C}}$ be the complexification of $G$.
Suppose there are given a Zariski connected $\mathbb{R}$-subgroup $K$ of $G$
with Lie algebra ${\mathfrak{k}}$ and connected unipotent $\mathbb{C}$
-subgroups $P_{\pm}$ of $G_{\mathbb{C}}$ with Lie algebras ${\mathfrak{p}}_{\pm}$. The group $G$ is called \textit{of Harish-Chandra type} if the
following conditions are satisfied:

(HC 1) $\mathfrak{g}_{\mathbb{C}}={\mathfrak{p}}_{+}+
{\mathfrak{k}}_{\mathbb{C}}\,+{\mathfrak{p}}_{-}$ is a direct sum of vector spaces, $\left[ {\mathfrak{k}}_{\mathbb{C}}{\mathfrak{,p}}_{\pm}\right] \subset{\mathfrak{p}}_{\pm}$, and $\overline{{\mathfrak{p}}_{+}}={\mathfrak{p}}_{-}$; (HC\ 2) the
map $P_{+}\times K_{\mathbb{C}}\times P_{-}$ $\rightarrow G_{\mathbb{C}}$ 
gives a holomorphic injection of $P_{+}\times K_{\mathbb{C}}\times P_{-}$
onto its open image $P_{+}K_{\mathbb{C}}P_{-}$; (HC\ 3) $G\subset P_{+}K_{\mathbb{C}}P_{-}$ and $G\cap K_{\mathbb{C}}P_{-}=K$.

If $g\in P_{+}K_{\mathbb{C}}P_{-}\subset G_{\mathbb{C}}$, we denote by
$(g)_{+}\in P_{+}$, $(g)_{0}\in K_{\mathbb{C}}$, $(g)_{-}\in P_{-}$ the
components of $g$ such that $g=(g)_{+}(g)_{0}(g)_{-}$. 

The identity connected component of a linear algebraic group $H$ is denoted
in the usual topology by $\mathring{H}$. The generalized Harish-Chandra
embedding \ of the homogeneous space $\mathcal{D}=\mathring{G}/\mathring{K}$
into ${\mathfrak{p}}_{+}$ is defined by $g\mathring{K}\longmapsto z$, where
 $g\in \mathring{G}$, $z\in $ ${\mathfrak{p}}_{+}$ and $\exp z=(g)_{+}$. 
Then the $\mathring{G}$-invariant complex structure of $\mathcal{D}$ is
determined by the natural inclusion $\mathcal{D}$ $\hookrightarrow P_{+}\subset
G_{\mathbb{C}}/(K_{\mathbb{C}}P_{-})$. Let $(G_{\mathbb{C}}\times 
{\mathfrak{p}}_{+})^{\prime }$ be the set of elements $(g,z)\in G_{\mathbb{C}}\times {\mathfrak{p}}_{+}$ such that $g\exp z\in P_{+}K_{\mathbb{C}}P_{-}$ and let $({\mathfrak{p}}_{+}\times {\mathfrak{p}}_{+})^{^{\prime }}$ be the set of
elements $(z_{1},z_{2})\in {\mathfrak{p}}_{+}\times {\mathfrak{p}}_{+}$ such
that $(\exp \bar{z}_2)^{-1}\exp z_{1}\in P_{+}K_{\mathbb{C}}P_{-}$.

The \emph{ canonical automorphy factor}
 $J:(G_{\mathbb{C}}\times{\mathfrak{p}}_{+})^{\prime}\rightarrow K_{\mathbb{C}}$
 and the \emph{canonical kernel function} $K:({\mathfrak{p}}_{+}\times{\mathfrak{p}}_{+})^{^{\prime}}\rightarrow K_{\mathbb{C}}$ for $G$ are defined by 
\begin{equation}
J(g,z)=\left( g\exp z\right) _{0},\quad K(z^{\prime},z)=
(((\exp\overline {z})^{-1}\exp z^{\prime})_0)^{-1},  \label{a4}
\end{equation}
where ($g,z)\in(G_{\mathbb{C}}\times{\mathfrak{p}}_{+})^{\prime}$ and 
$(z^{\prime},z)\in({\mathfrak{p}}_{+}\times{\mathfrak{p}}_{+})^{^{\prime}}$.

According to \cite{TA99}, \cite{LEE03} (Corollary 4.5, Proposition 4.7, and
equation (6.1)), we obtain

\begin{Theorem}\label{TH2.1}
a) The Jacobi group $G^{J}$  acts
transitively on $\mathfrak{D}^{J}$ by  
\begin{equation}
gx=(\sigma w,v_{\sigma w}+\,^{t}(c\tau(w)+d)^{-1}z),\ \rho(\sigma)=\left( 
\begin{array}{cc}
a & b \\ 
c & d
\end{array}
\right) ,  \label{a5}
\end{equation}
where $g=(\sigma,v,\varkappa)\in G^{J}$ and $x=(w,z)\in 
\mathfrak{D}^{J}$.

b) The canonical automorphy factor $J$ for the Jacobi group
 $G^{J}$ is given by 
\begin{equation}
J(g,x)=(J_{1}(\sigma ,w),0,J_{2}(g,x)),  \label{a6}
\end{equation}
where $g=(\sigma ,v,\varkappa )\in G^{J}$, $x=(w,z)\in 
\mathfrak{D}^{J}$, $J_{1}$ is the canonical automorphy factor 
for $G^{s}$, and 
\begin{equation}
J_{2}(g,x)=\varkappa +\frac{1}{2}D(v,v_{\sigma w})+\frac{1}{2}D(2v+\rho
(\sigma )z,J_{1}(\sigma ,w)z).  \label{aa6}
\end{equation}
c) The canonical kernel function $K$ for the Jacobi group 
$G^{J}$ is given by 
\begin{equation}
K(x,x^{\prime })=(K_{1}(w,w^{\prime }),0,K_{2}(x,x^{\prime })),  \label{a7}
\end{equation}%
where $x=(w,z)\in \mathfrak{D}^{J}$, $x^{\prime }=(w^{\prime
},z^{\prime })\in \mathfrak{D}^{J}$\emph{,} $K_{1}$ is the canonical
kernel function for $G^{s},$ and 
\begin{equation}
K_{2}(x,x^{\prime })=D(2\bar{z}^{\prime }+\frac{1}{2}\overline{\tau
(w^{\prime })}z,qz)+\frac{1}{2}D(\bar{z}^{\prime },q\tau (w)\bar{z}^{\prime
}),\ q=\rho (K_{1}(w,w^{\prime }))^{-1}.  \label{a8}
\end{equation}%
\end{Theorem}

The Heisenberg group $\mathrm{H}_{n}(\mathbb{R})$ consists of all elements
 $(\lambda ,\mu ,\kappa )$, where $\lambda ,\mu \in M_{1\,n}(\mathbb{R)}$, 
$\kappa \in \mathbb{R} $ with the multiplication law 
\begin{equation}
(\lambda ,\mu ,\kappa )\circ (\lambda ^{\prime },\mu ^{\prime },\kappa
^{\prime })=(\lambda +\lambda ^{\prime },\mu +\mu ^{\prime },\kappa +\kappa
^{\prime }\!+\lambda \,^{t}\!\mu ^{\prime }-\mu \,^{t}\!\lambda ^{\prime }).
\label{a10}
\end{equation}%
Let $G_{n}^{J}=\mathrm{Sp}(n,\mathbb{R})\ltimes \mathrm{H}_{n}(\mathbb{R})$
be the semidirect product of the symplectic group $\mathrm{Sp}(n,\mathbb{R})$
and the Heisenberg group $\mathrm{H}_{n}(\mathbb{R})$ endowed with the
following multiplication law: 
\begin{equation}
(\sigma ,(\lambda ,\mu ,\kappa ))\cdot (\sigma ^{\prime },(\lambda ^{\prime
},\mu ^{\prime },\kappa ^{\prime }))=\,(\sigma \sigma ^{\prime },(\lambda
\sigma ^{\prime },\mu \sigma ^{\prime },\kappa )\circ (\lambda ^{\prime
},\mu ^{\prime },\kappa ^{\prime })),  \label{a11}
\end{equation}%
where $(\lambda ,\mu ,\kappa )$,$\,(\lambda ^{\prime },\mu ^{\prime },\kappa
^{\prime })\in \mathrm{H}_{n}(\mathbb{R})$ and $\sigma ,\sigma ^{\prime }\in 
\mathrm{Sp}(n,\mathbb{R})$. The \emph{Jacobi group  $G_{n\,}^{J}$ of degree}
$n$ 
 acts transitively on the \emph{Jacobi-Siegel space}
$\mathfrak{H}_{n}^{J}=$ $\mathfrak{H}_{n}\times \mathbb{C}^{n}$ by $g(\Omega
,\zeta )=(\Omega _{g},\zeta _{g})$, where $(\Omega ,\zeta )\in \mathfrak{H}
_{n}^{J}$, $g=\left( \sigma ,\left( \lambda ,\,\mu , \kappa \right) \right)
\in G_{n}^{J}$, $\sigma $ is given by (\ref{a1}), and \cite{Y07} 
\begin{equation}
\Omega _{g}=(a\Omega +b)(c\Omega +d)^{-1},\ \zeta _{g}=\nu (c\Omega
+d)^{-1},\ \nu =\zeta +\lambda \Omega +\mu .  \label{b00}
\end{equation}

According with Theorem \ref{TH2.1} and \cite{SA71C}, we have

\begin{Proposition}\label{Pr2.2}
The canonical automorphy factor $J_{1}$ 
and the canonical kernel function $K_{1}$ for 
$\mathrm{Sp}(n,\mathbb{R}) $ are given by 
\begin{align}
J_{1}(\sigma,\Omega) & =\left( 
\begin{array}{cc}
\,^{t}(c\Omega+d)^{-1} & 0 \\ 
0 & c\Omega+d%
\end{array}
\right) ,\   \label{b1} \\
K_{1}(\Omega^{\prime},\Omega) & =\left( 
\begin{array}{cc}
0 & \overline{\Omega}-\Omega^{\prime} \\ 
(\Omega^{\prime}-\overline{\Omega})^{-1} & 0
\end{array}
\right) ,  \label{b2}
\end{align}
where $\Omega,\Omega^{\prime}\in\mathfrak{H}_{n}$ and 
$\sigma\in\mathrm{Sp}(n,\mathbb{R}) $ is given by {\emph{(\ref{a1})}}.

The canonical automorphy factor  $\theta=J_{2}(g,(\Omega,\zeta))$ for  $G^{J}_n$ is given by 
\begin{equation}
\theta=\kappa+\lambda\,\,{}^{t}\zeta+\nu\,\,{}^{t}\lambda-\nu(c\Omega
+d)^{-1}c\,{}^{t}\nu,\quad\nu=\zeta+\lambda\Omega+\mu,  \label{b4}
\end{equation}
where $g=\left( \sigma,\left( \lambda,\,\mu, \kappa\right) \right)
\in G_{n}^{J}$, $\sigma$ is given by {\emph{(\ref{a1})}}, and  $(\Omega
,\zeta)\in\mathfrak{H}_{n}^{J}$.

The canonical automorphy kernel $K_{2}$ \emph{for }$G^{J}_n$  is
given by  
\begin{equation}
K_{2}\left( (\zeta^{\prime},\Omega^{\prime}),(\zeta,\Omega)\right) =-\frac{1
}{2}(\zeta^{\prime}-\bar{\zeta})(\Omega^{\prime}-\bar{\Omega}
^{\prime})^{-1}(^{t}\zeta^{\prime}-\,^{t}\bar{\zeta}),  \label{b5}
\end{equation}
where $(\Omega,\zeta),(\Omega^{\prime},\zeta^{\prime})\in \mathfrak{H}
_{n}^{J}$.
\end{Proposition}

Let $\mathfrak{D}_{n}$ \emph{be the Siegel disk of degree} $n$ consisting of
all symmetric matrices $W \in M_{n}(\mathbb{C)}$ with $I_{n}-W\bar{W}>0$. Let $\mathrm{Sp}(n,\mathbb{R})_{\ast }$ be the multiplicative group of all
matrices $\omega \in M_{2n}(\mathbb{C)}$ such that
\begin{equation}
\omega =\left( 
\begin{matrix}
p & q \\ 
\bar{q} & \bar{p}
\end{matrix}
\right) ,\ \ ^{t}p\bar{p}-\,^{t}\bar{q}q=I_{n},\ ^{t}p\bar{q}=\,^{t}\bar{q}p
,\quad p,q\in M_{n}(\mathbb{C)}.  \label{jd1}
\end{equation}%
Remark that $\mathrm{Sp}(n,\mathbb{R})_{\ast }=\mathrm{Sp}(n,\mathbb{C})\cap 
\mathrm{U}(n,n)$ for $n>1$. $\mathrm{Sp}(n,\mathbb{R})_{\ast }$ acts
transitively on $\mathfrak{D}_{n}$ by $\omega W=(pW+q)(\bar{q}W+\bar{p}
)^{-1} $, where $\omega \in \mathrm{Sp}(n,\mathbb{R})_{\ast }$ and $\, W \in 
\mathfrak{D}_{n}$. \smallskip Let $K_{n\ast }\cong \mathrm{U}(n)$ be the
maximal compact subgroup of $\mathrm{Sp}(n,\mathbb{R})_{\ast }$ consisting
of all $\omega \in $ $\mathrm{Sp}(n,\mathbb{R})_{\ast }$ given by (\ref{jd1}
) with $p\in \mathrm{U}(n)$ and $q=0$. Then $\mathfrak{D}
_{n}\mathcal{\cong }\mathrm{Sp}(n,\mathbb{R})_{\ast }/\mathrm{U}(n)$.

Let $G_{n\ast }^{J}$ be the Jacobi group consisting of all elements $(\omega
,(\alpha ,\varkappa ))$, where $\omega \in $ $\mathrm{Sp}(n,\mathbb{R}
)_{\ast }$, $\alpha \in \mathbb{C}^{n}$, $\varkappa \in \mathrm{i}\mathbb{R}$,
 and endowed with the multiplication law 
\begin{equation}
(\omega ^{\prime },(\alpha ^{\prime },\varkappa ^{\prime }))(\omega ,(\alpha
,\varkappa ))=\left( \omega ^{\prime }\omega ,\varkappa +\varkappa ^{\prime
}+\beta {}^{t}\bar{\alpha}-\bar{\beta}{}^{t}\alpha \right) ,  \label{jd2}
\end{equation}%
where $(\omega ,(\alpha ,\varkappa )),(\omega ^{\prime },(\alpha ^{\prime
},\varkappa ^{\prime }))\in G_{n\ast }^{J}$ $\beta =$ $\alpha ^{\prime }p+
\bar{\alpha}^{\prime }{}\bar{q}$, and $\omega $ is given by (\ref{jd1}).

 The Heisenberg group 
$\mathrm{H}_{n}(\mathbb{R})_{\ast }$ consists of all elements $(I_{n},(\alpha ,\varkappa ))\in
G_{n\ast }^{J}$, where $\omega \in $ $\mathrm{Sp}(n,\mathbb{R})_{\ast }$, 
$\alpha \in \mathbb{C}^{n}$, $\varkappa \in \mathrm{i}\mathbb{R}$. The center 
$\mathcal{A}_{\ast }\mathcal{\cong }\mathbb{R}$ of $\mathrm{H}_{n}(\mathbb{R}
)_{\ast }$ consists of all elements $(I_{n},(0,\varkappa ))\in G_{n\ast
}^{J} $ with $\varkappa \in \mathrm{i}\mathbb{R}$. According to \cite{Y08}, there 
exists an isomorphism $\Theta :G_{n}^{J}\rightarrow $\ $G_{n\ast }^{J}$
given by $\Theta (g)=g_*$, $g=(\sigma ,(\lambda ,\mu ,\kappa ))\in G_{n}^{J}$, 
$g_* = (\omega ,(\alpha ,\varkappa ))\in G_{n\ast }^{J}$,
\begin{equation}
\sigma=\left( 
\begin{matrix}
a & b \\ 
c & d%
\end{matrix}
\right) , ~ \omega=\left( 
\begin{matrix}
p_{+} & p_{-} \\ 
\overline{p}_{-} & \overline{p}_{+}
\end{matrix}
\right) ,  \label{jd3}
\end{equation}%
\begin{equation}
p_{\pm}={\frac{1}{2}}(a\pm d)\pm\frac{\mathrm{\,i}}{2}\,(b\mp c),\quad
\alpha={\frac{1}{2}}(\lambda+ \mathrm{i}\mu),\quad\varkappa=-\mathrm{i}{\frac{\kappa}{2
}.}  \label{jd4}
\end{equation}

Let $\mathfrak{D}_{n}^{J}=\mathfrak{D}_{n}\times 
\mathbb{C}^{n}\mathcal{\cong }G_{n\ast }^{J}/(\mathrm{U}(n)\times \mathbb{R})$ be \emph{the Siegel-Jacobi disk} \emph{of degree n}. $G_{n\ast }^{J}$ acts transitively on 
$\mathfrak{D}_{n}^{J}$ by $g_{\ast }(W,z)=(W_{g_{\ast }},z_{g_{\ast }})$,
where $g_{\ast }=(\omega ,(\alpha ,\varkappa ))\in G_{n\ast}^{J}$, $(W,z)\in 
\mathfrak{D}_{n}^{J}$, $\omega $ is given by (\ref{jd3}), and \cite{Y08} 
\begin{equation}
W_{g_{\ast }}=(pW+q)(\bar{q}W+\bar{p})^{-1},~z_{g_{\ast }}=(z+\alpha W+
\bar{\alpha})(\bar{q}W+\bar{p})^{-1}.  \label{jd5}
\end{equation}

We now consider a partial Cayley transform of the Siegel-Jacobi disk $%
\mathfrak{D}_{n}^{J}$ onto the Siegel-Jacobi space $\mathfrak{H}_{n}^{J}$
which gives a partially bounded realization of $\mathfrak{H}_{n}^{J}$ 
\cite{Y08}.  The \emph{partial Cayley transform} $\phi:\mathfrak{D}%
_{n}^{J}\rightarrow\mathfrak{H}_{n}^{J}$ \ is defined by 
\begin{equation}
\Omega=\mathrm{i}(I_{n}+W)(I_{n}-W)^{-1},\quad\zeta=2\,
\mathrm{i}\,z(I_{n}-W)^{-1},  \label{jd6}
\end{equation}
where $(\zeta,\Omega)=\phi\left( (W,z)\right) $ and $(W,z)\in\mathfrak{D}
_{n}^{J}$.

$\phi$ is a a biholomorphic map which satisfies $g\phi=\phi g_{\ast}$ for
any $g\in G_{n}^{J}$ and $g_{\ast}=\Theta(g)$ \cite{Y08}.

The inverse partial Cayley transform $\phi^{-1}:\mathfrak{H}
_{n}^{J}\rightarrow\mathfrak{D}_{n}^{J}$ is given by 
\begin{equation}
W=(\Omega-\mathrm{i}I_{n})(\Omega+\mathrm{i}I_{n})^{-1},\quad z=\,\zeta
(\Omega+\mathrm{i}I_{n})^{-1},  \label{jd7}
\end{equation}
where $(W,z)=\phi^{-1}\left( (\Omega,\zeta)\right)\in \mathfrak{D}_{n}^{J} $ and $(\Omega,\zeta )\in\mathfrak{H}_{n}^{J}$.

According with Theorem \ref{TH2.1}, \cite{SA71C} and \cite{Y08}, we have

\begin{Proposition}\label{Pr2.3}
The canonical automorphy factor
$J_{1\ast}$ 
and the canonical kernel function $K_{1\ast}$ for $\mathrm{Sp}(n,
\mathbb{R})_{\ast}$  are given by
\begin{equation}
J_{1\ast}(\omega,W)=\left( 
\begin{array}{cc}
\,^{t}(\bar{q}W+\bar{p})^{-1} & 0 \\ 
0 & \bar{q}W+\bar{p}%
\end{array}
\right) ,  \label{jd8}
\end{equation}%
\begin{equation}
K_{1\ast}(W^{\prime},W)=\left( 
\begin{array}{cc}
I_{n}-W^{\prime}\bar{W} & 0 \\ 
0 & \,^{t}(I_{n}-W^{\prime}\bar{W})^{-1}
\end{array}
\right) ,  \label{jd9}
\end{equation}
where $W,W^{\prime}\in\mathfrak{D}_{n}$ and $\omega
 \in\mathrm{Sp}(n,\mathbb{R})_{\ast}$ is given by \emph{(\ref{jd1}).}

The canonical automorphy factor  $\theta _{\ast }=J_{2}(g_{\ast },(W,z))$
 for $G_{\ast }^{J}$ is given by
\begin{equation}
\theta _{\ast }=\kappa _{\ast }+z\,{}^{t}\alpha +\nu _{\ast }\,{}^{t}\alpha
-\nu _{\ast }(\bar{q}W+\bar{p})^{-1}\bar{q}\,{}^{t}\nu _{\ast },\quad
\nu _{\ast } = z+\alpha W+\bar{\alpha},\quad  \label{jd11}
\end{equation}%
where  $g_{\ast }=(\omega ,(\alpha ,\varkappa ))\in G_{n\ast }^{J}$, $
\omega $ is given by \emph{(\ref{jd1})}, and $(W,z)\in
 \mathfrak{D}_{n}^{J} $.

The canonical automorphy kernel for $G_{\ast }^{J}$  is given
by $K_{2\ast }((W^{\prime },z^{\prime }),(W,z))=A(W^{\prime },z^{\prime
};W,z)$, where $(W,z),(W^{\prime },z^{\prime })\in \mathfrak{D}
_{n}^{J} $, and 
\begin{equation}
A(W^{\prime },z^{\prime };W,z)=(\bar{z}+\frac{1}{2}z^{\prime }\bar{W}
)(I_{n}-W^{\prime }\bar{W})^{-1\,t}z^{\prime }+\frac{1}{2}\bar{z}
(I_{n}-W^{\prime }\bar{W})^{-1}W^{\prime }\,^{t}\bar{z}.  \label{jd12}
\end{equation}
\end{Proposition}

\section{Fock spaces based on the Siegel disk}\label{sec3}

Let $\mathcal{A}_{\ast }\mathcal{\cong }\mathbb{R}$ be the center of the
Heisenberg group $\mathrm{H}_{n\ast }(\mathbb{R})$. Given $m\in \mathbb{R}$,
let $\chi ^{m}$ be the central character of $\mathcal{A}_{\ast }$ defined by 
$\chi ^{m}(\kappa )=\exp \left( 2\pi \mathrm{i}m\kappa \right) $,$\ \kappa
\in \mathcal{A}_{\ast }$. Suppose $m>0$.

For each $W\in \mathfrak{D}_{n}$ we denote by $\mathcal{F}_{mW}$ the Fock
space of all functions $\Phi \in \mathcal{O}(\mathbb{C}^{n})$ such that $
\left\Vert \Phi \right\Vert _{mW}<\infty $ and the inner product is
defined by \cite{SA71C} 
\begin{eqnarray}
(\Phi ,\Psi )_{mW} &=&(2\pi m)^{n}\left( \det (1-W\bar{W})\right) ^{-1/2}
\label{h1} \\
&&\times \int_{\mathbb{C}^{n}}\Phi (z)\overline{\Psi (z)}\exp
(-8\pi mA(W,z))d\nu (z),\nonumber
\end{eqnarray}%
where the Lebesgue measure on $\mathbb{C}^{n}$ is given by 
\begin{equation}
d\nu (\zeta )=\pi ^{-n}\prod_{i=1}^{n}d 
 \mathrm{\operatorname{Re}}\zeta_{i} \,d
 \mathrm{\operatorname{Im}}\zeta _{i},  \label{h2}
\end{equation}%
and $A(W,z)=K_{2\ast }((W,z),(W,z))$ can be written as 
\begin{equation}
A(W,z)=(\bar{z}+\frac{1}{2}z\bar{W})(I_{n}-W\bar{W})^{-1\,t}z+\frac{1}{2}
\bar{z}(I_{n}-W\bar{W})^{-1}W\,^{t}\bar{z}.  \label{h3}
\end{equation}%

Remark that $\mathcal{F}_{m\,0}$ is the usual Bargmann space \cite{BAR61}.

  We
consider the Gaussian functions $G_{U}:\mathfrak{D}_{n}^{J}\rightarrow 
\mathbb{C}$, $U\in \mathbb{C}^{n}$, defined by $G_{U}(W,Z)=G(U,Z,W)$, where

\begin{equation}
G(U,Z,W)=\exp (U\,^{t}Z+\frac{1}{2}UW\,^{t}U)=\sum_{s\in \mathbb{N}^{n}}
\frac{U^{s}}{s!}P_{s}(Z,W)  \label{h5}
\end{equation}%
for all $(W,Z)\in \mathcal{D}_{n}^{J}$. We utilize the notation 
\begin{equation}
U^{s}={\prod_{i=1}^{n}U_{i}^{s_{i}}},\quad s!={\prod_{i=1}^{n}}
s_{i}!,\quad |s|={\sum_{i=1}^{n}}s_{i},\quad \delta _{sr}=
{\prod_{i=1}^{n}}\delta _{s_{i}r_{i}},  \label{h6}
\end{equation}%
where $U=(U_{1},...,U_{n})\in M_{1n}(\mathbb{C)\cong C}^{n}$, 
$s=(s_{1},...,s_{n})\in \mathbb{N}^{n}$, and $r=(r_{1},...,r_{n})\in 
\mathbb{N}^{n}$.\ The polynomials $P_{s}:\mathcal{D}_{n}^{J}\rightarrow \mathbb{C}$, 
$s\in \mathbb{N}^{n}$, are exactly the matching functions studied by Neretin 
\cite{NE90}. We express the homogeneous polynomial $P_{s}$ of degree $|s|$
in the following compact form: 
\begin{equation}
P_{s}(Z,W)=\sum_{a\in A_{n},\,\tilde{a}\leq s}
\frac{s!}{2^{\hat{a}}a!(s-\tilde{a})!}Z^{s-\tilde{a}}W^{a},  \label{h7}
\end{equation}%
where $A_{n}$ is the set of all symmetric matrices $a=(a_{ij})_{1\leq
i,\,j\leq n}$ with $a_{ij}\in \mathbb{N}$, 
\begin{equation}
W^{a}={\prod_{1\leq i\leq j\leq n}W_{ij}^{\,a_{ij}}},\quad a!=
{\prod_{1\leq i\leq j\leq n}}a_{ij},\quad \tilde{a}_{k}=
\sum_{i=1}^{n}a_{ik},\quad \hat{a}=\sum_{i=1}^{n}a_{ii},  \label{h8}
\end{equation}%
and $\,\tilde{a}\leq s$ is equivalent with $\,\tilde{a}_{i}\leq s_{i}$ for 
$1\leq i\leq n$. Using the equations

\begin{equation}
\int_{\mathbb{C}^{n}}U^{s}\bar{U}^{r}d\nu (U)=\delta _{sr}s!,
\label{h10}
\end{equation}%
\begin{equation}
\!\int_{\mathbb{C}^{n}}\!\!G(U,Z^{\prime },W^{\prime })G(\bar{U},\bar{Z},\bar{W
})d\nu (U\!)\!=\!\det (1\!-\!W^{\prime }\!\bar{W})^{-1/2}\!\!\exp\! A(W^{\prime }\!,z^{\prime
};W\!,\!z),  \label{h9}
\end{equation}
where $A(W^{\prime },z^{\prime };W,z)$ is defined by (\ref{jd12}), we obtain 
\begin{equation}
\left( \det (1-W^{\prime }\bar{W})\right) ^{-1/2}\exp A(W^{\prime}
\!,z^{\prime }\!;W,z)=\sum_{s\in \mathbb{N}^{n}}\frac{1}{s!}P_{s}(Z^{\prime
},W^{\prime })\overline{P_{s}(Z,W)}.  \label{h11}
\end{equation}

Equation (\ref{h9}) is given in \cite{I67} (Lemma 5).

We now define the polynomials $\Phi _{W\!s}\in \mathcal{F}_{mW}$, $s\in 
\mathbb{N}^{n}$, by 
\begin{equation}
\Phi _{W\!s}(z)=\frac{1}{\sqrt{s!}}P_{s}(2\sqrt{2\pi m}z,W),\ s\in 
\mathbb{N}^{n}.  \label{h26}
\end{equation}

\begin{Proposition}\label{Pr3.1}
Given  $W\in \mathfrak{D}_{n}$,  the set
of polynomials $\{\Phi _{W\!s}|s\in \mathbb{N}^{n}\}$ forms an
orthonormal basis of the Fock space $\mathfrak{F}_{mW}.$  The kernel
function of  $\mathfrak{F}_{mW}$  admits the expansion 
\begin{equation}
\left( \det (1-W\bar{W})\right) ^{-1/2}\exp \left( 2\pi mA(W,z^{\prime
};W,z)\right) =\sum_{s\in \mathbb{N}^{n}}\Phi _{W\!s}(z^{\prime })
\overline{\Phi _{W\!s}(z)}.  \label{h27}
\end{equation}
\end{Proposition}

\emph{Proof}. Given $U\in\mathbb{C}^{n}$ and $W\in\mathcal{D}_{n}$, we
define the function $\Psi_{UW\,}:\mathbb{C}^{n}\rightarrow\mathbb{C}$ such
that $\Psi_{UW\,}(z)=G(U,2\sqrt{2\pi m}z,W)$. Using the change of variables 
$Z=$ $2\sqrt{2\pi m}z$, we have
\begin{equation}
\!\left\Vert\! \Psi_{UW\,}\right\Vert _{mW}^{2}\!=\!\pi^{-n}\!\det(1\!-\!W\bar{W}
)^{-1/2}\!\int_{\mathbb{C}^{n}}\!\exp(B(U,Z,W)\!-\!A(Z,W))d\nu(Z)\,,
\label{h28}
\end{equation}
where
\begin{equation}
B(U,Z,W)=U^{t}Z+\bar{U}Z^{\dag }-\frac{1}{2}UW\,^{t}U-\frac{1}{2}\bar{U}
\bar{W}U^{+}.  \label{h29}
\end{equation}
Using the change of variables $Z^{\prime }=(1-W\bar{W})^{-1/2}
\left( Z-\bar{U}-WU\right) $, the relation $d\nu (Z)=\det (1-W\bar{W})d\nu (Z^{\prime })$,
and the relation \cite{BAR61} 
\begin{equation}
\int_{\mathbb{C}^{n}}\!\!\exp\!(\! -\!{\overline{Z^{\prime}}}
^{t}\!Z'\!-\!\frac{1}{2}(Z'\bar{W}^{t}Z'\!+\!\overline{
Z}^{\prime }WZ^{\prime \dag }))\! d\nu (Z^{\prime })\!=\!\pi ^{n}\!\left(\! \det
(1\!-\!W\bar{W})\right)\! ^{-1/2},  \label{h30}
\end{equation}%
we obtain $\left\Vert \Psi _{UW\,}\right\Vert _{mW}^{2}=\exp (UU^{\dag })$.
Then 
\begin{equation}
\sum_{s,r\in \mathbb{N}^{n}}\frac{U^{s}\bar{U}^{r}}{\sqrt{s!r!}}(\Phi
_{W\!s},\Phi _{W\!r})_{mW}=\sum_{s\in \mathbb{N}^{n}}\frac{1}{s!}U^{s}
\bar{U}^{s}.  \label{h31}
\end{equation}%
By comparing the coefficients of $U^{s}\bar{U}^{r}$ in the series of both
sides of (\ref{h31}), we see that 
\begin{equation}
(\Phi _{W\!s},\Phi _{W\!r})_{mW}=\delta _{sr}s!,\quad s,r\in \mathbb{N}
^{n}.  \label{h32}
\end{equation}%
Using (\ref{h11}) and (\ref{h26}), we obtain the expansion (\ref{h27}). \hfill $\gata$

We now introduce the set of polynomials $f_{s}:\mathcal{D}_{n}^{J}\rightarrow
\mathbb{C}$, $s\in\mathbb{N}^{n}$, defined by 
\begin{equation}
f_{s}(W,z)=\frac{1}{\sqrt{s!}}P_{s}(2\sqrt{2\pi m}z,W).  \label{h33}
\end{equation}

Let $\mathcal{H}_{0}(\mathcal{D}_{n}^{J})$ be the complex linear subspace of
all holomorphic functions $f\in\mathcal{O}(\mathcal{D}_{n}^{J})$ with the
basis $\{f_{s}|s\in\mathbb{N}^{n}\}$. Let $\mathfrak{F}_{m}(\mathcal{D}
_{n}^{J})$ be the Hilbert space of all functions $f\in\mathcal{O}(\mathcal{D}
_{n}^{J})$ such that $\left\langle f,f\right\rangle _{m\,}<\infty$, where
the inner product $\left\langle .,.\right\rangle _{m\,}$is defined such that
the set $\{f_{s}|s\in\mathbb{N}^{n}\}$ is an orthonormal basis. We now prove

\begin{Proposition}\label{Pr3.2}
 a) The generating function of the basis
$\{f_{s}|s\in\mathbb{N}^{n}\}$  can be expressed as
\begin{equation}
\exp(8\pi mU\,^{t}z+\frac{1}{2}UW\,^{t}U)=\sum_{s\in\mathbb{N}^{n}}
\frac{U^{s}}{\sqrt{s!}}f_{s}(W,z).  \label{h34}
\end{equation}

The kernel function of $\mathfrak{F}_{m}(\mathcal{D}_{n}^{J}) $ 
admits the expansion 
\begin{equation}
\left( \det (1-W^{\prime }\bar{W})\right) ^{-1/2}\exp A(W^{\prime
},z^{\prime };W,z)=\sum_{s\in \mathbb{N}^{n}}f_{s}(W^{\prime },z^{\prime })
\overline{f_{s}(W,z)}.  \label{h35}
\end{equation}

b) $f$ $\in\mathcal{O}(\mathcal{D}_{n}^{J})$ is a solution of
the system of differential equations
\begin{equation}
\frac{\partial^{2}f}{\partial z_{j}\partial z_{k}}=8\pi m(1+\delta_{jk})
\frac{\partial f}{\partial W_{jk}},\ 1\leq j\leq k\leq n,  \label{h36}
\end{equation}
if and only if $f$ $\in\mathcal{H}_{0}(\mathcal{D}_{n}^{J})$.
\end{Proposition}

\emph{Proof. }Using\emph{\ }(\ref{h5}) and (\ref{h33}), we obtain
 (\ref{h34}). The generating function (\ref{h34}) satisfies (\ref{h36}). Then 
\begin{equation}
\frac{\partial^{2}f_{s}}{\partial z_{j}\partial z_{k}}=8\pi m(1+\delta _{jk})
\frac{\partial f_{s}}{\partial W_{jk}},\ 1\leq j\leq k\leq n,\
 s\in\mathbb{N}^{n}.  \label{h37}
\end{equation}
\emph{\ } Using (\ref{h7}) and (\ref{h33}), we obtain
\begin{equation}
f_{s}(z,W)=\frac{1}{\sqrt{s!}}(2\sqrt{2\pi m}z)^{s}+R_{s}(z,W),  \label{h38}
\end{equation}
where $R_{s}$ is a polynomial of degree $|s|-1$ in $z$. Then there exists
the change of basis $\{z^{s}W^{a}|s\in\mathbb{N},a\in A_{n}\}\longmapsto
\{f_{s}(z,W)W^{a}|s\in\mathbb{N}^{n},a\in A_{n}\}$ in 
$\mathcal{O}(\mathcal{D}_{n}^{J})$. Let $f$ $\in\mathcal{O}(\mathcal{D}_{n}^{J})$. Then there
exists the set $\{c_{s}|c_{s}\in\mathcal{O}(\mathcal{D}_{n}),\ s\in\mathbb{N}
^{n}\}$ such that
\begin{equation}
f(z,W)=\sum_{s\in\mathbb{N}^{n}}c_{s}(W)f_{s}(W,z).  \label{h39}
\end{equation}
If $f$ satisfies (\ref{h36}), then $\partial c_{s}/\partial W_{jk}=0$ for
any $1\leq j\leq k\leq n,\ s\in\mathbb{N}^{n}$. Then $c_{s}$ is constant for
any $s\in\mathbb{N}^{n}$. Hence $f$ $\in\mathcal{O}(\mathcal{D}_{n}^{J})$.
The inverse implication follows from (\ref{h37}). \hfill  $\gata$

In the case $n=1$, Proposition \ref{Pr3.2} has been obtained in \cite{BERC08B}.

\section{Scalar holomorphic discrete series of the \\ Jacobi group on the Siegel-Jacobi disk}\label{sec4}

Consider the Jacobi group $G_{n}^{J}$.
Let $\delta$ be a rational representation of $GL(n,\mathbb{C})$ such that
$\delta|_{\text{U}(n)}$ is a scalar irreducible representation of the unitary
group U$(n)$ with highest weight $k$, $k\in\mathbb{Z}$, and $\delta(A)=(\det A)^{k}$ \cite{Z89}. 
 Let $m\in \mathbb{R}$. Let $\chi
=\delta \otimes {\bar{\chi}}^{m}$, where the central character $\chi ^{m}$
of $\mathcal{A\cong }\mathbb{R}$ is defined by $\chi ^{m}(\kappa )=\exp
\left( 2\pi \mathrm{i}m\kappa \right) $,$\ \kappa \in \mathcal{A}$.  Any scalar holomorphic irreducible
representation of $G_{n}^{J}$ is characterized by an index $m$ and a weight 
$k$. Suppose $m>0$ and $k>n+1/2$.

Let $\mathcal{H}^{mk}$ denote the Hilbert space of all holomorphic
functions $\varphi\in\mathcal{O}(\mathfrak{H}_{n}^{J}\mathfrak{)}$ such that 
$\left\Vert \varphi\right\Vert _{\mathfrak{H}_{n}^{J}}<\infty$ with the
inner product defined by \cite{TA90}

\begin{equation}
(\varphi ,\psi )_{\mathfrak{H}_{n}^{J}}=C\int_{\mathfrak{H}_{n}^{J}}\varphi
(\Omega ,\zeta )\overline{\psi (\Omega ,\zeta )}\,\mathcal{K}^{mk}(\Omega
,\zeta )^{-1}d\mu (\Omega ,\zeta ),\   \label{4.2}
\end{equation}%
where $C$ is a positive constant, $(\Omega ,\zeta )\in \mathfrak{H}_{n}^{J}$
and the $G_{n}^{J}$-invariant measure on $\mathfrak{H}_{n}^{J}$ is given by 
\begin{equation}
d\mu (\Omega ,\zeta )=(\det Y)^{-n-2}\prod_{1\leq i\leq n}d\xi _{i}\,d\eta
_{i}\prod_{1\leq j\leq k\leq n}dX_{jk}\,dY_{jk}. \ \   \label{4.3}
\end{equation}
Here $ \xi =\operatorname{Re} \zeta,~ \eta = \operatorname{Im} \zeta, ~ X = 
\operatorname{Re} \Omega,~ Y = \operatorname{Im} \Omega$.

The kernel function $\mathcal{K}^{mk}$ is defined by \cite{TA90} 
\begin{equation}
\mathcal{K}^{mk}(\Omega ,\zeta)=\mathsf{K}^{mk}((\Omega ,\zeta
),(\Omega ,\zeta ))=\exp \left( 4\pi m\eta Y^{-1\ t}\eta \right) (\det
Y)^{k},  \label{4.4}
\end{equation}%
\begin{equation}
\mathsf{K}^{mk}((\zeta ^{\prime },\Omega ^{\prime }),(\zeta ,\Omega
))\!=\!\left(\! \det (\frac{\mathrm{i}}{2}\bar{\Omega}-\frac{\mathrm{i}}{2}\Omega
^{\prime })\!\right)\! ^{-k}\!\exp\! \left( 2\pi \mathrm{im}K\left( (\zeta ^{\prime
},\Omega ^{\prime }),(\zeta ,\Omega )\right) \right) ,  \label{4.4b}
\end{equation}
where $K$ is given by (\ref{b5}).

Let $\pi^{mk}$ be the unitary representation of $G_{n}^{J}$ on $\mathcal{H}
^{mk}$ defined by \cite{TA90} 
\begin{equation}
\left( \pi^{mk}(g^{-1})\varphi\right) (\Omega,\zeta)=\mathcal{J}
^{mk}(g,(\Omega,\zeta))\varphi(\Omega_{g},\zeta_{g}),  \label{4.5}
\end{equation}
where $\varphi\in\mathcal{H}^{mk}$, $g\in G_{n}^{J}$, $(\Omega,\zeta )\in
\mathfrak{H}_{n}^{J}$ and $(\Omega_{g},\zeta_{g})\in\mathfrak{H}_{n}^{J}$ is
given by (\ref{b00}).

The automorphic factor $\mathcal{J}^{mk}$ for $G_{n}^{J}$ is defined by 
\cite{TA90} 
\begin{equation}
\mathcal{J}^{mk}(g,(\zeta ,\Omega ))=\left( \det (c\Omega +d)\right)
^{-k}\exp (2\pi \mathrm{i}m\theta ),  \label{4.7}
\end{equation}%
where $\theta $ is given by (\ref{b4}) and $\sigma $ is given by  
(\ref{a1}).

Takase proved the following theorem \cite{TA90}, \cite{TA92}:

\begin{Theorem}\label{TH4.1}
Suppose $k>n+1/2$. Then $\mathcal{H}
^{mk}\neq\{0\}$ and $\pi^{mk}$ is an irreducible unitary
representation of $G_{n}^{J}$ which is square integrable modulo center.
\end{Theorem}

Let $\mathcal{H}_{\ast }^{mk}$ denote the complex pre-Hilbert space of all 
$\psi \in \mathcal{O}(\mathfrak{D}_{n}^{J})$ such that $
\left\Vert \psi \right\Vert _{\mathfrak{D}_{n}^{J}}<\infty $ with the
inner product defined by 
\begin{equation}
(\psi _{1},\psi _{2})_{\mathfrak{D}_{n}^{J}}=C_{\ast }\int_{\mathfrak{D}
_{n}^{J}}\psi _{1}(W,z)\,\overline{\psi _{2}(W,z)}\left( \mathcal{K}_{\ast
}^{mk}(W,z)\right) ^{-1}d\,\nu (W,z),  \label{4.9}
\end{equation}
where $C_{\ast }$ is a positive constant, $(z,W)\in \mathfrak{D}_{n}^{J}$, 
\begin{equation}
\mathcal{K}_{\ast }^{mk}(W,z)=\left( \det (I_{n}-W\bar{W})\right)
^{-k}\exp (8\pi mA(W,z)),  \label{4.10}
\end{equation}%
where $A$ is given by (\ref{h3}) and the $G_{n}^{J}$-invariant measure on 
$\mathfrak{D}_{n}^{J}$ is \cite{Y08}
\begin{equation}
d\nu (W\!,z)\!=\!(\det (1\!-\!W\bar{W}))^{-n-2}\!\prod_{i=1}^{n}\!d\,\mathrm{\operatorname{
Re}}z_{i}\,d\,\mathrm{\operatorname{Im}}z_{i}\!\prod_{1\leq j\leq k\leq n}\!\!d 
\operatorname{Re}W_{jk}d\,\mathrm{\operatorname{Im}}W_{jk.}  \label{4.12}
\end{equation}

According with \cite{SA71B}, \cite{Y08}, and (\ref{jd12}), the kernel
function $\mathcal{K}_{\ast }^{mk}$ is given by $\mathcal{K}_{\ast
}^{mk}(W,z)=\mathsf{K}_{\ast }^{mk}((W,z),(W,z))$, where 
\begin{equation}
\mathsf{K}_{\ast }^{mk}((z,W),(z^{\prime },W^{\prime }))\!=\!\left( \det
(I_{n}-W^{\prime }\bar{W})\right) ^{-k}\!\exp \left( 8\pi mA(W^{\prime
},z^{\prime };W,z)\right) .  \label{4.12b}
\end{equation}

\begin{Remark}\label{Rem4.2}
{\emph { Using the coherent state method, Kramer, Saraceno, and
Berceanu investigated the kernel (\ref{4.10}) in the case $8\pi m=1$ 
\cite{KRSAR82}, \cite{BERC06}-\cite{BERC08B}.}}
\end{Remark}

We now introduce the map $g_{\ast}\longmapsto\pi_{\ast}^{mk}(g_{\ast})$,
where $\pi_{\ast}^{mk}(g_{\ast})$: $\mathcal{H}_{\ast}^{mk}\rightarrow$ 
$\mathcal{H}_{\ast}^{mk}$ is defined by 
\begin{equation}
\left( \pi_{\ast}^{mk}(g_{\ast}^{-1})\psi\right)
(z,W)=J_{\ast}^{mk}(g_{\ast},(z,W))\psi(z_{g_{\ast}},W_{g_{\ast}}),
\label{4.13}
\end{equation}
$\psi\in\mathcal{H}_{\ast}^{mk},\ g_{\ast}=(\omega,(\alpha,\varkappa))\in
G_{n\ast}^{J}$ , $(z,W)\in\mathfrak{D}_{n}^{J}$, and $(z_{g_{\ast}},W_{g_{%
\ast}})\in\mathfrak{D}_{n}^{J}$ is given by (\ref{jd5}). The automorphic
factor $J_{\ast}^{mk}$ for $G_{n\ast}^{J}$ is defined by \cite{SA71B}, 
\cite{Y08} 
\begin{equation}
J_{\ast}^{mk}(g_{\ast},(z,W))=\exp(2\pi\mathrm{i}m\theta_{\ast})\left(
\det(\bar{q}W+\bar{p})\right) ^{-k}  \label{4.15}
\end{equation}
where $\theta_{\ast}$ is given by (\ref{jd11}) and $\omega$ given by 
(\ref{jd1}).

\begin{Proposition}\label{Pr4.3}
Suppose  $m>0$, $k>n+1/2$, and 
$C=2^{n(n+3)}C_{\ast}$. Then

a) $\mathcal{H}_{\ast}^{mk}\neq\{0\}$ and $\pi_{\ast}^{mk}$
is an irreducible unitary representation of $G_{n\ast}^{J}$ on
the Hilbert space $\mathcal{H}_{\ast}^{mk}$ which is square
integrable modulo center.

b) There exists the unitary isomorphism $T_{\ast}^{mk}:\mathcal{H}_{\ast}^{mk}\rightarrow\mathcal{H}^{mk}$ given by 
\begin{equation}
\,\varphi(\Omega,\zeta)=\psi\,\left( W,z\right) (\det(I_{n}-W))^{k}\exp(4\pi
mz(I_{n}-W)^{-1}\,^{t}z),  \label{4.17}
\end{equation}
where $\psi\in\mathcal{H}_{\ast}^{mk}$, \ $\varphi=T_{\ast
}^{mk}(\psi)$, $\left( W,z\right) \in\mathfrak{D}_{n}^{J}$, $
(\Omega,\zeta)=\phi(\left( -W,z\right) )\in\mathfrak{H}_{n}^{J}$, 
and $\phi$ is given by \emph{(\ref{jd6})}.

The inverse isomorphism $T^{mk}:\mathcal{H}^{mk}\rightarrow 
\mathcal{H}_{\ast}^{mk}$  is given by
\begin{equation}
\psi\,\left( W,z\right) =\varphi(\Omega,\zeta)\,\,(\det(I_{n}-\mathrm{i}%
\Omega))^{k}\exp\left( 2\pi m\zeta(I_{n}-\mathrm{i}\Omega)\,^{-1}\,^{t}\zeta
\right) \,,\   \label{4.18}
\end{equation}
where $\psi\in\mathcal{H}^{mk}_{\ast}$\emph{,} \ $\psi=T^{mk}(\varphi )$
\emph{,} $(\Omega,\zeta)\in\mathfrak{H}_{n}^{J}$\emph{, }$\left( -W,z\right)
=\phi^{-1}\left( (\Omega,\zeta)\right) \in\mathfrak{D}_{n}^{J}$, and
 $\phi^{-1}$ is given by \emph{(\ref{jd7})}.

c) The representations $\pi^{mk}$ and  $\pi_{\ast}^{mk}$ 
are unitarily equivalent. 
\end{Proposition}

\emph{Proof. }Using the partial Cayley transform (\ref{jd6}) and 
(\ref{jd7}), we obtain 
\begin{equation}
Y=\operatorname{Im}\Omega=(I_{n}-W)^{-1}(I_{n}-W\bar{W})(I_{n}-\bar{W})^{-1},
\label{4.19}
\end{equation}%
\begin{equation}
\eta=\operatorname{Im}\zeta=z(I_{n}-W)^{-1}+\bar{z}(I_{n}-\bar{W})^{-1}.
\label{4.20}
\end{equation}
By (\ref{4.19}) and (\ref{4.20}), we obtain
\begin{equation}
\eta Y^{-1\ t}\eta=2A(z,-W)-z(I_{n}-W)^{-1}\,^{t}z-\bar{z}(I_{n}-\bar{W}
)^{-1}\,^{t}\bar{z},  \label{4.21}
\end{equation}
where $A$ is given by (\ref{h3}). Using (\ref{4.3}), (\ref{4.12}),
 (\ref{jd6}), and (\ref{jd7}), in the limit $\Omega\rightarrow\mathrm{i}I_{n}$ and $W\rightarrow0$, we obtain 
\begin{equation}
d\mu(\zeta,\Omega)=2^{n(n+3)}d\nu(z,W).  \label{4.22}
\end{equation}

By (\ref{4.2}), (\ref{4.9}), (\ref{4.17}), (\ref{4.18}), the condition
 $C=2^{n(n+3)}C_{\ast}$, and the change of variables $W\rightarrow-W$, we get 
$\left\Vert \varphi\right\Vert _{\mathfrak{H}_{n}^{J}}=\left\Vert
\psi\right\Vert _{\mathfrak{D}_{n}^{J}}$. From $\zeta(I_{n}-\mathrm{i}
\Omega)\,^{-1}\,^{t}\zeta=-2\,z(I_{n}-W)^{-1}\,^{t}z$ it is clear that 
(\ref{4.17}) and (\ref{4.18}) are equivalent. By Theorem 4.1, a) and b) hold.
Using (\ref{jd6}), (\ref{jd7}), (\ref{4.5}), (\ref{4.13}), (\ref{4.17}), and
(\ref{4.18}), we obtain $\pi^{mk}T^{mk}=T_{\ast}^{mk}\pi_{\ast}^{mk}$. \hfill $\gata$

\begin{Remark}\label{Re4.4}
{\emph{Berndt,  B\"{o}cherer  and Schmidt constructed the holomorphic
discrete series of the Jacobi group in the case $n=1$ \cite{BB90}, \cite{BS98}.}}
\end{Remark}

Let $\mathcal{H}^{\,k}$ denote the complex Hilbert space of all holomorphic
functions $\Phi \in \mathcal{O}(\mathfrak{D}_{n}\mathfrak{)}$ such that $\left\Vert
\Phi \right\Vert _{\mathfrak{D}_{n}}<\infty $, with the inner product
defined by 
\begin{align}
(\Psi _{1},\Psi _{2})_{k}& =\int_{\mathcal{D}_{n}}\Psi _{1}(W)\overline{\Psi
_{2}(W})\left( \det (1-W\bar{W})\right) ^{k-1/2}d\mu _{\mathcal{D}_{n}}(W),
\label{4.24} \\
d\mu _{\mathcal{D}_{n}}(W)& =\left( \det (1-W\bar{W})\right)
^{-n-1}\prod_{1\leq j\leq k\leq n}\mathrm{d}\operatorname{Re}W_{jk}\,\mathrm{d}\operatorname{Im}W_{jk}.  \notag
\end{align}
$\ $We have $\mathcal{H}^{k}\neq \{0\}$ for $k>n+1/2$ \cite{BE75}, \cite{TA90}. Let $\{Q_{a}|a\in A_{n}\}$ be an orthonormal polynomial basis of $\mathcal{H}^{\,k}$.

We introduce the polynomials 
\begin{equation}
F_{sa}(W,z)=\sqrt{\frac{\left( 8\pi m\right) ^{n}}{C_{\ast }s!}}
P_{s}(\sqrt{8\pi m}z,W)Q_{a}(W),\ s\in \mathbb{N}^{n},\ a\in A_{n}.  \label{4.25}
\end{equation}
\begin{Proposition}\label{Pr4.5}
The set of polynomials $\left\{
F_{sa}|s\in \mathbb{N}^{n},\ a\in A_{n}\right\} $ forms an
orthonormal basis of $\mathcal{H}_{\ast }^{mk}$. The kernel
function of $\mathcal{H}_{\ast }^{mk}$ satisfies the expansion 
\begin{equation}
( \det (1-W^{\prime }\bar{W}))\!^{-k}\exp A(W^{\prime },z^{\prime
},W,z)\!=\!\sum_{s\in \mathbb{N}^{n}, a\in A_{n}}\!F_{sa}(W^{\prime },z^{\prime })
\overline{F_{sa}(W,z)}.  \label{4.26}
\end{equation}
\end{Proposition}

\emph{Proof. }We introduce the functions $F_{U\,}:\mathfrak{D}
_{n}^{J}\rightarrow \mathbb{C}$, $U\in \mathbb{C}^{n}$, such that $
F_{U\,}(W,z)=G(U,2\sqrt{2\pi m}z,W)$. Using (\ref{4.25}) and the proof of
Proposition 3.1, we have \emph{\ } 
\begin{align}
(F_{U\,}(W,z)Q_{a}\!,\!F_{U\,}(W,z)Q_{b})_{\mathfrak{D}_{n}^{J}}& \!=C_{\ast
}\left( 8m\right) ^{n}\exp (UU^{\dag })  \label{4.27} \\
& \!\times \!\int_{\mathcal{D}_{n}}\!Q_{a}(W)\!\overline{Q_{b}(W)}\!\det (\!1\!-\!W\bar{W}%
)^{k-1/2}\!d\mu _{\mathcal{D}_{n}}(W),  \notag
\end{align}%
\begin{equation}
(F_{sa}\,,F_{rb})_{\mathfrak{D}_{n}^{J}}=\delta _{sr}\delta _{ab},\
s,r\in \mathbb{N}^{n},\ a,b\in A_{n}.  \label{4.28}
\end{equation}%
The Berezin kernel of $\mathcal{H}^{\,k}$ is positive definite for $k>n+1/2$ 
\cite{BE75} and satisfies the following identity: \emph{\ } 
\begin{equation}
\left( \det (1-W^{\prime }\bar{W})\right) ^{-k+1/2}=\sum_{\ a\in
A_{n}}Q_{a}(W^{\prime })\overline{Q_{a}(W)}.  \label{4.29}
\end{equation}%
Using (\ref{h27}) and (\ref{4.29}), we obtain (\ref{4.26}). \hfill $\gata$

\begin{Remark}\label{Re4.6}
 {\emph{In the case $n=1$ and $8\pi m=1$, the
expansion (\ref{4.26}) was obtained in \cite{BERC06B}, using the coherent state method. }}
\end{Remark}

\begin{Remark}\label{Re4.7}
{\emph{ We now discuss the unitary representations of Jacobi
groups based on Siegel-Jacobi domains in the language of coherent states 
\cite{P86}. Let $Q(\mathcal{H)}$ be the set of all one-dimensional
projections of the Hilbert space $\mathcal{H}$. Let $P[\psi ]$ denote the
one-dimensional projection determined by $\psi \in \mathcal{H}\backslash
\{0\}$. The elements of $Q(\mathcal{H)}$ can be considered either as normal
pure states of the von Neumann algebra of bounded operators on $\mathcal{H}$
or as pure states of the $C^{\ast }$-algebra of compact operators on 
$\mathcal{H}$ \cite{C83}. The projective  Hilbert space $P(\mathcal{H})$ consists of
all one-dimensional complex linear subspaces of $\mathcal{H}$. The space 
$P(\mathcal{H})$ is a K\"{a}hler manifold equipped with the usual Fubini-Study
metric \cite{C83}. The space $Q(\mathcal{H)}$ with relative $w^{\ast }$
-topology is homeomorphic to $P(\mathcal{H})$ with the manifold topology 
\cite{C83}. Then we can identify $Q(\mathcal{H)}$ with $P(\mathcal{H})$.}}

{\emph{We recall an intrinsic definition of coherent state representations given in  \cite{M78}.}}

{\emph{Let $G$ be a connected, simply connected Lie group and $\mathfrak{X}$ a $G$-homogeneous space which admits an invariant measure $\mu_{\mathfrak{X}}$. Let $\pi$ be
a continuous irreducible unitary representation of $G$ in the separable
Hilbert space $\mathcal{H}$. A family $\mathcal{E}=\left\{ E_{x}|x\in
\mathfrak{X}\right\} $ of one-dimensional projections in $\mathcal{H}$ will
be called a $\pi$-\emph{system of coherent states} based on $\mathfrak{X}$
if the following conditions are satisfied:\ 1) $E_{g\,x}=\pi\,(g)E_{x}\pi%
\,(g)^{-1}$ for any $g\in G$ and $x\in\mathfrak{X}$; 2) there exists $\psi\in
\mathcal{H}\backslash\{0\}$, such that $\int _{\mathfrak{X}}\left\vert \left\langle
\psi,\pi\,(g)\psi\right\rangle \right\vert ^{2}d\,\mu_{\mathfrak{X}}<\infty$. $\pi$ is
called a \emph{symplectic (K\"{a}hler) coherent state representation} if $
\mathcal{E}$ and $\mathfrak{X}$ are isomorphic symplectic (K\"{a}hler)
manifolds and $\mathfrak{X}$ is a symplectic (K\"{a}hler) submanifold of 
$Q(\mathcal{H})$.}}

{\emph{Moscovici and Verona have been studied coherent state representations based
precisely on the coadjoint orbit associated with $\pi$ in the sense of
geometric quantization \cite{M78}. The Schr\"{o}dinger coherent state
systems for the Heisenberg group with one-dimensional center on the Fock
spaces of holomorphic functions have been obtained by Bargmann \cite{BAR61},
Satake \cite{SA71}, \cite{SA71C}, \cite{SA80}, and Lee \cite{LEE03}.}}

{\emph{Lisiecki and Neeb investigated some K\"{a}hler coherent state
representations of Heisenberg groups and Jacobi groups with one-dimensional
center \cite{LIS91},\cite{N00}. The orbit method for the Heisenberg group
and the Jacobi group with multi-dimensional center has been studied in
detail by Yang \cite{Y02}.}}
  
{\emph{Let $\pi$ be an irreducible unitary representation of the Jacobi group 
$G^{J} $ with the Jacobi-Siegel domain $\mathfrak{D}$ and the kernel function 
$K:\mathfrak{D\times D\rightarrow}\mathrm{Hom}(W,W)$. The representation
space $\mathcal{H}$ consists of holomorphic functions taking their values in
a finite dimensional Hilbert space $W$. For each $x\in\mathfrak{D}$ and 
$v\in W$, we consider the vectors $K_{xv}\in\mathcal{H}$ given by $
K_{xv}(x^{\prime })=K(x,x^{\prime})v$ for any $x^{\prime}\in$ $\mathfrak{D}$. Then $\{P[K_{xv}]|x\in\mathfrak{D}, ~v\in W\}$ is a $\pi$-system of coherent
states. In particular, the $\pi^{mk}$-system of coherent states based on 
$\mathfrak{H}_{n}^{J}$ and the $\pi_{\ast}^{mk}$-system of coherent states
based on $\mathfrak{D}_{n}^{J}$ are determined by the explicit kernel functions given by (\ref{4.4b}) and (\ref{4.12b}), respectively.}}
\end{Remark}
\textbf{Acknowledgments.} {\small  The authors are indebted to the Organizing Committee of {\it The International Conference
          of Differential Geometry and Dynamical Systems}, University Politehnica of Bucharest, Romania, 
             for the opportunity to report results at the meeting. S. B. was partially supported by
the CNCSIS-UEFISCSU project PNII- IDEI 454/2009, CNCSIS Cod  ID-44, A.G. was
partially supported by the CNCSIS-UEFISCSU project PNII- IDEI 545/2009, CNCSIS Cod
 ID-1089 and both the authors have been partially supported by the
ANCS project  program PN 09 37 01 02/2009.}

\end{document}